\documentclass[a4paper,12pt]{article}

\usepackage{amsmath}
\usepackage{amssymb}
\usepackage{amsfonts}
\usepackage{graphicx}

\title{\bf Optimal switching control design for polynomial systems:
an LMI approach}

\begin{document}

\author{Didier Henrion$^{1,2,3}$, Jamal Daafouz$^4$, Mathieu Claeys$^1$}

\footnotetext[1]{CNRS; LAAS; 7 avenue du colonel Roche, F-31077 Toulouse; France. {\tt henrion@laas.fr}}
\footnotetext[2]{Universit\'e de Toulouse; UPS, INSA, INP, ISAE; UT1, UTM, LAAS; F-31077 Toulouse; France}
\footnotetext[3]{Faculty of Electrical Engineering, Czech Technical University in Prague,
Technick\'a 2, CZ-16626 Prague, Czech Republic}
\footnotetext[4]{Universit\'e de Lorraine, CRAN, CNRS, IUF, 2 avenue de la forêt de Haye, 54516 Vand\oe uvre cedex, France. {\tt jamal.daafouz@univ-lorraine.fr}}

\date{Draft of \today}

\maketitle

\begin{abstract}
We propose a new LMI approach to the design of optimal switching sequences
for polynomial dynamical systems with state constraints. We formulate
the switching design problem as an optimal control problem which is
then relaxed to a linear programming (LP) problem in the space
of occupation measures. This infinite-dimensional LP can be solved numerically
and approximately with a hierarchy of convex finite-dimensional LMIs.
In contrast with most of the existing work on LMI methods, we have
a guarantee of global optimality, in the sense that we obtain an asympotically converging
(i.e. with vanishing conservatism) hierarchy of lower bounds on the achievable performance.
We also explain how to construct an almost optimal switching sequence.
\end{abstract}

\section{Introduction}
A switched system is a particular class of a hybrid system that consists of a set of dynamical subsystems, one of which is active at any instant of time, and a policy for activating and deactivating the subsystems. One may encounter such dynamical systems in a wide variety of application domains such as automotive industry, power systems, aircraft and traffic control, and more generally the area of embedded systems. Switched systems have been the concern of many researchers and many results are available for stability analysis and control design. They put in evidence the important fact that it is possible to orchestrate the subsystems through an adequate switching strategy in order to impose global stability. Interested readers may refer to the survey papers \cite{DeCarlo, Liberzonb, Shorten, Lin} and the interesting and useful books \cite{Liberzon, Sun} and the references therein. 

In this context, switching plays a major role for stability and performance properties. Indeed, switched systems are generally controlled by switched controllers and the control signal is intrinsically discontinuous. As far as optimality is concerned, several results are also available in two main different contexts:
\begin{itemize}
\item the first category of methods exploits necessary optimality conditions, in the form
of Pontryagin's maximum principle (the so-called indirect approaches), or through a large nonlinear
discretization of the problem (the so-called direct approaches), see
\cite{Bengea, Branicky, Cassandras, Hed99, Piccoli, Riedinger99, Riedinger, Seatzu, Sussmann, Shaikh2007, Xu2} for details. Therefore only local optimality can be guaranteed for general
nonlinear problems, even when discretization can be properly controlled;
\item the second category collects extensions of the performance indexes $H_2$ and $H_\infty$
originally developped for linear time invariant systems without switching, and use the flexibility
of Lyapunov's approach, see for instance \cite{Geromel, Deaecto} and references therein. Even for linear switched systems, the proposed results are based on nonconvex optimization problems (e.g. bilinear matrix inequality conditions) difficult to solve directly. Sufficient linear matrix inequality (LMI) design
conditions may be obtained, but at the price of introducing a conservatism (pessimism) which is
hard, if not impossible, to evaluate.
Since the computation of this optimal strategy is a difficult task, a suboptimal solution is of interest  only when it is proved to be consistent, meaning that it imposes to the switched system a performance
not worse than the one produced by each isolated subsystem \cite{Gero}. 
\end{itemize}

Despite the interest of these existing approaches, the optimal control problem is not completely solved for switched systems and new strategies are more than welcome, as computationally viable design techniques are missing. In this paper, we consider the problem of designing optimal switching rules in the case of polynomial switched dynamical systems. Classically, we formulate the optimal control switching problem as an optimal control problem with controls being functions of time valued in $\{0,1\}$,
and we relax it into a control problem with controls being functions of time valued in $[0,1]$.
In contrast with existing approaches following this relaxation strategy, see e.g. \cite{Bengea, Riedinger99},
relying on Pontryagin's maximum principle,
our aim is to apply the approach of \cite{sicon}, which consists
in modeling control and trajectory functions as occupation measures.
This allows for a convex linear programming (LP) formulation of the optimal control problem.
This infinite-dimensional LP can be solved numerically and approximately with a hierarchy of convex finite-dimensional LMIs. On the one hand, our approach follows the optimal control modeling framework.
On the other hand, it exploits the flexibility and computational efficiency of
the convex LMI framework. In contrast with most of the existing work on LMI methods, we have
a guarantee of global optimality, in the sense that we obtain an asympotically converging
(i.e. with vanishing conservatism) hierarchy of lower bounds on the achievable performance.

The paper is organized as follows: in Section \ref{pb} we state the optimal switching problem
to be solved, and we propose an alternative, relaxed formulation allowing chattering
of the trajectories. In Section \ref{meas}
we introduce occupation measures as a device to linearize the optimal control problem
into an LP problem in the cone of nonnegative measures.
In Section \ref{solvelp} we explain how to solve the resulting infinite-dimensional LP
with a converging hierarchy of finite-dimensional LMI problems. As explained
in Section \ref{seq}, an approximate optimal switching sequence can be extracted
from the solutions of the LMI problem, and this is illustrated with classical
examples in Section \ref{exs}. Finally, the paper ends with a sketch
of further research lines.

\section{Optimal switching problem}\label{pb}

Consider the optimal control problem
\begin{equation}\label{ocp}
\begin{array}{rcll}
p^* & = & \inf & \int_0^T l_{\sigma(t)}(t,x(t))dt \\
&& \mathrm{s.t.} & \dot{x}(t) = f_{\sigma(t)}(t,x(t)),\quad \sigma(t) \in \{1,2,\ldots,m\}\\
&&& x(0) \in X_0, \quad x(T) \in X_T\\
&&& x(t) \in X, \quad t \in [0,T]
\end{array}
\end{equation}
with given polynomial velocity field $f_{\sigma(t)} \in {\mathbb R}[t,x]^n$
and given polynomial Lagrangian $l_{\sigma(t)} \in {\mathbb R}[t,x]$ indexed
by an integer-valued signal $\sigma : [0,T]\to\mathbb \{1,2,\ldots,m\}$.
System state $x(t)$ belongs to a given compact semialgebraic set $X \subset {\mathbb R}^n$
for all $t \in [0,T]$
and the initial state $x(0)$ resp. final state $x(T)$ are constrained to a given
compact semialgebraic set $X_0 \subset X$ resp. $X_T \subset X$.
In problem (\ref{ocp}) the infimum is w.r.t. sequence $\sigma$ and terminal time $T$.

In this paper, for the sake of simplicity,
we assume that the terminal time $T$ is finite, that is, we do not
consider the asymptotic behavior. Typically, if a solution
to problem (\ref{ocp}) is expected for a very large or infinite terminal time,
we must reformulate the problem by relaxing the state constraints.

Optimal control problem (\ref{ocp}) can then be equivalently written as
\begin{equation}\label{discrete}
\begin{array}{rcll}
p^* & = & \inf_{u} & \int_0^T \sum_{k=1}^m l_k(t,x(t))u_k(t)dt \\
&& \mathrm{s.t.} & dx(t) = \sum_{k=1}^m f_k(t,x(t))u_k(t)dt\\
&&& x(0) \in X_0, \quad x(T) \in X_T\\
&&& x(t) \in X, \quad u(t) \in U, \quad t \in [0,T]
\end{array}
\end{equation}
where the infimum is with respect to a time-varying vector $u(t)$
which belongs for all $t \in [0,T]$ to the (nonconvex) discrete set
\[
U := \{(1,0,\ldots,0), (0,1,\ldots,0), \ldots, (0,0,\ldots,1)\} \subset {\mathbb R}^m.
\]
In general the infimum in problem (\ref{discrete}) is not attained (see our
numerical examples later on) and the problem
is relaxed to
\begin{equation}\label{convex}
\begin{array}{rcll}
p^*_R & = & \inf & \int_0^T \sum_{k=1}^m l_k(t,x(t))u_k(t)dt \\
&& \mathrm{s.t.} & dx(t) = \sum_{k=1}^m f_k(t,x(t))u_k(t)dt\\
&&& x(0) \in X_0, \quad x(T) \in X_T\\
&&& x(t) \in X, \quad u(t) \in \mathrm{conv}\:U, \quad t \in [0,T]
\end{array}
\end{equation}
where the minimization is now with respect to a time-varying vector
$u(t)$ which belongs for all $t \in [0,T]$ to the (convex) simplex
\[
\mathrm{conv}\:U = \{u \in {\mathbb R}^m \: :\: \sum_{k=1}^m u_k = 1, \quad u_k \geq 0, \quad k=1,\ldots,m\}.
\]
In \cite{Bengea}, problem (\ref{convex}) is called the embedding of problem (\ref{ocp}),
and it is proved that the set of trajectories
of problem (\ref{ocp}) is dense (w.r.t. the uniform norm in the space of continuous
functions of time) in the set of trajectories
of embedded problem (\ref{convex}). Note however that these authors consider
the more general problem of switching design in the presence
of additional bounded controls in each individual dynamics.
To cope with chattering effects due to the simultaneous presence
of controls and (initial and terminal) state constraints,
they have to introduce a further relaxation of the embedded
control problem. In this paper, we do not have controls in
the dynamics, and the only design parameter is the switching
sequence.

An equivalent way of writing the dynamics in problem (\ref{convex})
is via a differential inclusion
\begin{equation}\label{inclusion}
\dot{x}(t) \in \mathrm{conv}\:\{f_1(t,x(t)),\ldots,f_m(t,x(t))\}.
\end{equation}
By this it is meant that at time $t$, the state velocity $\dot{x}(t)$
can be any convex combination of the vector fields $f_k(t,x(t))$,
$k=1,\ldots,m$, see e.g. \cite[Section 3.1]{bressan} for
a tutorial introduction. For this reason, problem (\ref{convex})
is also sometimes called the convexification
of problem (\ref{ocp}). 

Since problem (\ref{convex}) is a relaxation of problem (\ref{ocp}),
it holds $p^*_R \leq p^*$. For most of the physically relevant problems,
and especially when the state constraints in problem (\ref{ocp})
are not overly stringent, it actually holds that $p^*_R=p^*$.
For a discussion about the cases for which $p^*_R<p^*$,
please refer to \cite[Appendix C]{roa} and references therein.

\section{Occupation measures}\label{meas}

Given an initial condition $x_0 \in X_0$
and an admissible control $u(t)$, denote by $x(t|x_0,u)$, $t \in [0,T]$,
the corresponding admissible trajectory, an absolutely continuous
function of time with values in $X$.
Define the occupation measure
\[
\mu(A\times B|x_0,u) := \int_0^T I_{A\times B}(t,x(t|x_0,u))dt
\]
for all subsets $A\times B$ in the Borel $\sigma$-algebra of subsets of $[0,T]\times X$,
where $I_A(x)$ is the indicator function of set $A$, equal to one if $x \in A$, and zero otherwise. 
We write $x(t|x_0,u)$ resp. $\mu(dt,dx|x_0,u)$ to emphasize the dependence of $x$ resp. $\mu$
on initial condition $x_0$ and control $u$, but for conciseness we also use the notation
$x(t)$ resp. $\mu(dt,dx)$. The occupation measure can be disintegrated into
\[
\mu(A\times B) = \int_A \xi(B|t)\omega(dt)
\]
where $\xi(dx|t)$ is the distribution of $x \in {\mathbb R}^n$, conditional on $t$,
and $\omega(dt)$ is the marginal w.r.t. time $t$, which models the control action
as a measure on $[0,T]$.
The conditional $\xi$ is a stochastic kernel, in the sense that for all $t \in [0,T]$,
$\xi(.|t)$ is a probability measure on $X$, and for every $B$ in the Borel $\sigma$-algebra
of subsets of $X$, $\xi(B|.)$ is a Borel measurable function on $[0,T]$.
An equivalent definition is $\xi(B|t) = I_B(x(t)) = \delta_{x(t)}(B)$ where $\delta$
is the Dirac measure. The occupation measure encodes the system trajectory, and
the value $\int_0^T \mu(dt,B) = \mu([0,T]\times B)$ is equal to the total time spent
by the trajectory in set $B \subset X$. Note also that time integration of
any smooth test function $v : [0,T]\times X \to {\mathbb R}$ along a trajectory
becomes a time and space integration against $\mu$, i.e.
\[
\int_0^T v(t,x(t))dt = \int_0^T \int_X v(t,x)\mu(dt,dx) = \int v \mu.
\]

In optimal control problem (\ref{convex}), we associate an occupation measure
\[
\mu_k(dt,dx)=\xi_k(dx|t)\omega_k(t)
\]
for each system mode $k=1,\ldots,m$, so that globally
\[
\sum_{k=1}^m \mu_k = \mu
\]
is the occupation measure of a system trajectory subject to switching.
The marginal $\omega_k$ is the control, modeled as a measure which is
absolutely continuous w.r.t. the Lebesgue measure, i.e. such that
\[
\int_X \mu_k(dt,dx) = \omega_k(dt) = u_k(t)dt
\]
for some measurable control function $u_k(t)$, $k=1,\ldots,m$.
The system dynamics in problem (\ref{convex}) can then be expressed as
\begin{equation}\label{omega}
dx(t) = \sum_{k=1}^m f_k(x(t)) \omega_k(dt)
\end{equation}
where the controls are now measures $\omega_k$.
To enforce that $u(t) \in \mathrm{conv}\:U$ for almost all times $t \in [0,T]$,
we add the constraint
\begin{equation}\label{lebesgue}
\sum_{k=1}^m \omega_k(dt) = I_{[0,T]}(t)dt
\end{equation}
where the right hand side is the Lebesgue measure, or uniform measure on $[0,T]$.

Given a smooth test function $v : [0,T]\times X \to {\mathbb R}$ and an admissible trajectory
$x(t)$ with occupation measure $\mu(dt,dx)$, it holds
\[
\begin{array}{l}
\int_0^T dv(t,x(t)) = v(T,x(T)) - v(0,x(0))  \\
\:\: =  \int_0^T \left(\frac{\partial v}{\partial t}(t,x(t))+
\sum_k \mathrm{grad}\:v(t,x(t)) f_k(t,x(t)) u_k(t) \right)dt\\
\:\: =\int_0^T \int_X \left(\frac{\partial v}{\partial t}(t,x)\mu(dt,dx)+\right.\\
\quad\quad \left.\sum_k \mathrm{grad}\:v(t,x) f_k(t,x) \mu_k(dt,dx) \right) \\
\:\: = \sum_k \int \frac{\partial v}{\partial t}\mu_k + \mathrm{grad}\:v f_k \mu_k.
\end{array}
\]
Now, consider that the initial state is not a single vector $x_0$ but a random vector
whose distribution is ruled by a probability measure
$\mu_0$, so that at time $t$ the state $x(t)$ is also modeled by a probability measure
$\mu_t(.) := \xi(.|t)$, not necessarily equal to $\delta_{x(t)}$.
The interpretation is that $\mu_t(B)$ is the probability that the state $x(t)$
belongs to a set $B \subset X$.
Optimal control problem (\ref{convex}) can then be formulated
as a linear programming (LP) problem:
\begin{equation}\label{lp}
\begin{array}{l}
p^*_M \:=\: \inf \:\sum_k \int l_k \mu_k \:\:\mathrm{s.t.} \\
\quad \int v \mu_T - \int v \mu_0 = 
\sum_k \int \frac{\partial v}{\partial t} \mu_k +\mathrm{grad}\:v  f_k \mu_k\\
\quad \sum_k \int w \mu_k = \int w\\
\quad \forall v \in C^1([0,T]\times X), \:\:\forall w \in C^1([0,T]) \\

\end{array}
\end{equation}
where the infimum is w.r.t. measures $\mu_0 \in M_+(X_0)$, $\mu_T \in M_+(X_T)$,
$\mu_k \in M_+([0,T]\times X)$, $k=1,\ldots,m$
with $M_+(A)$ denoting the cone of finite nonnegative measures supported on $A$,
identified as the dual of the cone of nonnegative continuous functions supported on $A$.

It follows readily that $p^*\geq p^*_R \geq p^*_M$, and under some
additional assumptions it should be possible to prove that $p^*_R=p^*_M$
and that the marginal densities $u_k$ extracted from solutions $\mu_k$
of problem (\ref{lp}) are optimal for problem (\ref{discrete})
and hence problem (\ref{ocp}). We leave the rigorous statement and
its proof for an extended version of this paper.

Note that the use of relaxations and LP formulations of optimal control
problems (on ordinary differential equations and partial differential equations)
is classical, and can be traced back to the work by L. C. Young, Filippov, and
then Warga and Gamkrelidze, amongst many others. For more details and a historical
survey, see e.g. \cite[Part III]{Fattorini}.

\section{Solving the LP on measures}\label{solvelp}

To summarize, we have formulated our relaxed optimal switching control
problem (\ref{convex}) as the convex LP (\ref{lp}) in the space of measures.
This can be seen as an extension of the approach of \cite{sicon} which was originally
designed for classical optimal control problems. Alternatively, this
can also be understood as an application of the approach of \cite{impulse}
where the control measures are restricted to be absolutely continuous w.r.t. time.
Indeed, absolute continuity of the control measures is enforced by relation
(\ref{lebesgue}).
The infinite-dimensional LP on measures (\ref{lp}) can be solved approximately by
a hierarchy of finite-dimensional linear matrix inequality (LMI) problems,
see \cite{sicon,impulse,roa} for details
(not reproduced here). The main idea behind the hierarchy is to manipulate
each measure via its moments truncated to degree $2d$, where $d$ is a relaxation
order, and to use necessary LMI conditions for a vector to contain moments
of a measure. The hierarchy then consists of LMI problems of increasing sizes,
and it provides a sequence of lower bounds
$p^*_d \leq p^*$ which is monotonically increasing, i.e. $p^*_d \leq p^*_{d+1}$
and asymptotically converging, i.e. $\lim_{d\to\infty} p^*_d = p^*$.

The number of variables $N_d$ at the LMI relaxation of order $d$
grows linearly in $m$ (the number of modes), and polynomially in $d$,
but the exponent is a linear function of $n$ (the number of states).
In practice, given the current state-of-the-art in general-purpose
LMI solvers and personal computers, we can expect an LMI problem
to be solved in a matter of a few minutes provided the problem
is reasonably well-conditioned and $N_d \leq 5000$.

\section{Optimal switching sequence}\label{seq}

Let
\[
y_{k,\alpha} := \int_0^T \int_X t^\alpha \mu_k(dt,dx) = \int_0^T t^\alpha \omega_k(dt),
\:\: \alpha=0,1,\ldots
\]
denote the moments of measure $\omega_k$, $k=1,\ldots,m+1$.
Solving the LMI relaxation of order $d$ yields real numbers
$\{y^{d}_{k,\alpha}\}_{\alpha=0,1,\ldots,2d}$ which are approximations to $y_{k,\alpha}$.

In particular, the zero order moment of each measure $\mu_k$
is an approximation of its mass, and hence of the time $t_k = \int \mu_k \in [0,T]$
spent by an optimal switching sequence on mode $k$. At each LMI relaxation $d$, it holds
$\sum_{k=1}^{m+1} y^{d}_{k,0} = 1$ for all $d$, and $\lim_{d\to\infty} y^{d}_{k,0} = t_k$,
so that in practice, it is expected that good approximations of $t_k$ are obtained at
relatively small relaxation orders.

The (approximate) higher order moments $\{y^d_{k,\alpha}\}_{\alpha=1,\ldots,2d}$
allow to recover (approximately) the densities $u_k(t)$ of each measure $\omega_k(dt)$,
for $k=1,\ldots,m+1$.
The problem of recovering a density from its moments is a well-studied
inverse problem of numerical analysis. Since we expect in many cases
the density to be piecewise constant, with possible discontinuities here corresponding to commutations
between system modes, we propose the following strategy.

Let us assume that we have the moments
\[
y_{\alpha} := \int_0^T t^{\alpha} \omega(dt)
\]
of a (nonnegative) measure with piecewise constant density
\[
\omega(dt) = u(t)dt := \sum_{k=1}^N u_k I_{[t_{k-1},t_k]}(t) dt
\]
such that the boundary values are zero, i.e. $u_0=0$ and $u_{N+1}=0$.
The Radon-Nikodym derivative of this measure reads
\[
u'(dt) = \sum_{k=1}^N (u_{k+1}-u_k) \delta_{t_k}(dt)
\]
where $\delta_{t_k}$ denotes the Dirac measure at $t=t_k$.
Let
\[
y'_{\alpha} := \int_0^T t^{\alpha} u'(dt) = \sum_{k=1}^N (u_{k+1}-u_k) t^{\alpha}_k
\]
denote the moments of the (signed) derivative measure. By integration by parts it holds
\[
y'_{\alpha} = -\alpha y_{\alpha-1}, \quad \alpha = 0,1,2,\ldots
\]
which shows that the moments of $u'$ can be obtained readily from the moments of $u$.
Since $u'$ is a sum of $N$ Dirac measures, the moment matrix of $u'$ is a (signed)
sum of $N$ rank-one moment matrices, and the atoms $t_k$ as well as the weights
$u_{k+1}-u_k$, $k=1,\ldots,N$ can be obtained readily from an eigenvalue
decomposition of the moment matrix as explained e.g. in \cite[Section 4.3]{lasserre}.

More generally, the reader interested in numerical methods for reconstructing
a measure from the knowledge of its moments is referred to
\cite{mevissen} and references therein, as well as to the recent works
\cite{decastro,batenkov,candes} which deal with the problem of reconstructing
a piecewise-smooth function from its Fourier coefficients.
To make the connection between moments and Fourier coefficients,
let us just mention that the moments $y_\alpha = \int t^\alpha u(t)dt$
of a smooth density $u(t)$ are (up to scaling) the Taylor coefficients
of the Fourier transform $\hat{u}(s) := \int e^{-2\pi i s t} u(t)dt =
\sum_\alpha \frac{(-2\pi i)^\alpha}{\alpha!} y_\alpha s^\alpha$.
If the $y_\alpha$ are given, then $\hat{u}(s)$ is given by its Taylor series,
and the density $u(t)$ is recovered with
the inverse Fourier transform $u(t) = \int e^{2\pi i s t} \hat{u}(s)ds$.
Numerically, an approximate density can be obtained by applying the inverse fast Fourier
transform to the (suitably scaled) sequence $\{y^d_\alpha\}_{\alpha=0,1,\ldots,2d}$ of moments.

\section{Examples}\label{exs}

\subsection{First example}

Consider the scalar $(n=1)$ optimal control problem (\ref{ocp}):
\[
\begin{array}{rcll}
p^* & = & \inf & \int_0^1 x^2(t)dt \\
&& \mathrm{s.t.} & \dot{x}(t) = a_{\sigma(t)}x(t)\\
&&& x(0) = \frac{1}{2}, \quad x(1) \in [-1,1]\\
&&& x(t) \in [-1,1], \quad\forall t \in [0,1]
\end{array}
\]
where the infimum is w.r.t. to a switching sequence
$\sigma : [0,1] \mapsto \{1,\:2\}$ and
\[
a_1 := -1, \quad a_2 := 1.
\]

In Table \ref{ex1table} we report the lower bounds $p^*_d$ on
the optimal value $p^*$ obtained by solving LMI relaxations of
increasing orders $d$, rounded to 5 significant digits.
We also indicate the number of variables
(i.e. total number of moments) of each LMI problem, as well as
the zeroth order moment of each occupation measure
(recall that these are approximations of the time spent
on each mode). We observe that the values of the
lower bounds and the masses stabilize quickly.

\begin{table}[h!]
\centering
\begin{tabular}{c|c|c|cc}
$d$ & $p^*_d$ & $N_d$ & $y^d_{1,0}$ & $y^d_{2,0}$ \\ \hline
1 & $-5.9672\cdot10^{-9}$ & 18 & 0.74056 & 0.25944 \\
2 & $4.1001\cdot10^{-2}$ & 45 & 0.75170 & 0.24830 \\
3 & $4.1649\cdot10^{-2}$ & 84 & 0.74632 & 0.25368 \\
4 & $4.1666\cdot10^{-2}$ & 135 & 0.74918 & 0.25082 \\
5 & $4.1667\cdot10^{-2}$ & 198 & 0.74974 & 0.25026 \\
6 & $4.1667\cdot10^{-2}$ & 273 & 0.74990 & 0.25010 \\
7 & $4.1667\cdot10^{-2}$ & 360 & 0.74996 & 0.25004 \\
\end{tabular}
\caption{Lower bounds $p^*_d$ on the optimal value $p^*$ obtained
by solving LMI relaxations of increasing orders $d$;
$N_d$ is the number of variables in the LMI problem;
$y^d_{k,0}$ is the approximate time spent on each mode $k=1,2$.\label{ex1table}}
\end{table}

In this simple case, it is easy to obtain analytically the
optimal switching sequence: it consists of driving the state
from $x(0)=\frac{1}{2}$ to $x(\frac{1}{2})=0$ with the first
mode, i.e. $u_1(t)=1,u_2(t)=0$ for $t \in [0,\frac{1}{2}[$, and
then chattering between the first and second mode with
equal proportion so as to keep $x(t)=0$, i.e. $u_1(t)=\frac{1}{2},u_2(t)=\frac{1}{2}$
for $t \in ]\frac{1}{2},1]$. It follows that the infimum
is equal to
\[
p^* = \int_0^{1/2} \left(\frac{1}{2}-t\right)^2 dt = \frac{1}{24} \approx 4.1667\cdot10^{-2}.
\]
Because of chattering, the infimum in problem (\ref{ocp})
is not attained by an admissible
switching sequence. It is however
attained in the convexified problem (\ref{convex}).

The optimal moments can be obtained analytically
\[
y_{1,{\alpha}} = \int_0^{\frac{1}{2}} t^{\alpha} dt + \frac{1}{2} \int_{\frac{1}{2}}^1
t^{\alpha} dt = \frac{2+2^{-\alpha}}{4+4\alpha},
\]
\[
y_{2,{\alpha}} = \frac{1}{2} \int_{\frac{1}{2}}^1 t^{\alpha} dt
= \frac{2-2^{-\alpha}}{4+4\alpha}
\]
and they can be compared with the following moment vectors
obtained numerically at the 7th LMI relaxation:
\[
\begin{array}{rcl}
y^7_1 & = & \left[0.74996 \:\: 0.31246 \:\: 0.18746 \:\: 0.13277 \:\: 0.10308 \cdots \right], \\
y_1 & = & \left[0.75000 \:\: 0.31250 \:\: 0.18750 \:\: 0.13281 \:\: 0.10313 \cdots \right], \\
y^7_2 & = & \left[0.25004 \:\: 0.18754 \:\: 0.14588 \:\: 0.11723 \:\: 0.096919 \cdots \right], \\
y_2 & = & \left[0.25000 \:\: 0.18750 \:\: 0.14583 \:\: 0.11719 \:\: 0.096875 \cdots \right].
\end{array}
\]
We observe that the approximate moments $y^7_k$ closely
match the optimal moments $y_k$, so that the approximate
control law $u_k$ extracted from $y^7_k$ will be
almost optimal.

\subsection{Second example}

We revisit the double integrator example with state constraint
studied in \cite{sicon}, formulated as the following
optimal switching problem:
\[
\begin{array}{rcll}
p^* & = & \inf & T \\
&& \mathrm{s.t.} & \dot{x}(t) = f_{\sigma(t)}(x(t))\\
&&& x(0) = [1,1], \quad x(T) = [0,0]\\
&&& x_2(t) \geq -1, \quad\forall t \in [0,T]
\end{array}
\]
where the infimum is w.r.t. to a switching sequence
$\sigma : [0,T] \mapsto \{1,\:2\}$ with free terminal time $T\geq 0$
and affine dynamics
\[
f_1 := \left[\begin{array}{c}x_2\\-1\end{array}\right], \quad
f_2 := \left[\begin{array}{c}x_2\\1\end{array}\right].
\]
We know from \cite{sicon} that the optimal sequence consists of starting with
mode $1$, i.e. $u_1(t)=1$, $u_2(t)=0$ for $t \in [0,2]$,
then chattering with equal proportion between mode $1$ and $2$,
i.e. $u_1(t)=u_2(t)=\frac{1}{2}$ for $t \in [2,\frac{5}{2}]$
and then eventually driving the state to the origin
with mode $2$, i.e. $u_1(t)=0$, $u_2(t)=1$ for $t \in [\frac{5}{2},\frac{7}{2}]$.
Here too the infimum $p^*=\frac{7}{2}$ is not attained for
problem (\ref{ocp}), whereas it is attained with the
above controls for problem (\ref{convex}).

In Table \ref{ex1table} we report the lower bounds $p^*_d$ on
the optimal value $p^*$ obtained by solving LMI relaxations of
increasing orders $d$, rounded to 5 significant digits.
We also indicate the number of variables
(i.e. total number of moments) of each LMI problem, as well as
the zeroth order moment of each occupation measure
(recall that these are approximations of the time spent
on each mode). We observe that the values of the
lower bounds and the masses stabilize quickly
to the optimal values $p^*=\frac{7}{2}$, $y_{1,0}=\frac{5}{2}$,
$y_{2,0}=\frac{5}{4}$.

\begin{table}[h!]
\centering
\begin{tabular}{c|c|c|cc}
$d$ & $p^*_d$ & $N_d$ & $y^d_{1,0}$ & $y^d_{2,0}$ \\ \hline
1 & $2.5000$ & 30 & 1.7500 & 0.75000 \\
2 & $3.2015$ & 105 & 2.1008 & 1.1008 \\
3 & $3.4876$ & 252 & 2.2438 & 1.2438 \\
4 & $3.4967$ & 495 & 2.2484 & 1.2484 \\
5 & $3.4988$ & 858 & 2.2494 & 1.2494 \\
6 & $3.4993$ & 1365 & 2.2496 & 1.2497 \\
7 & $3.4996$ & 2040 & 2.2498 & 1.2498 \\
\end{tabular}
\caption{Lower bounds $p^*_d$ on the optimal value $p^*$ obtained
by solving LMI relaxations of increasing orders $d$;
$N_d$ is the number of variables in the LMI problem;
$y^d_{k,0}$ is the approximate time spent on each mode $k=1,2$.\label{ex2table}}
\end{table}

The optimal switching sequence, which corresponds here to
control measures $\omega_k(dt)$ with piecewise constant densities,
is obtained numerically as explained in Section \ref{seq}, by considering the moments
of the (weak) derivative of control measures.

\subsection{Third example}

Consider the optimal control problem (\ref{ocp}):
\[
\begin{array}{rcll}
p^* & = & \inf & \int_0^{\infty} \|x(t)\|^2_2 dt \\
&& \mathrm{s.t.} & \dot{x}(t) = A_{\sigma(t)}x(t)\\
&&& x(0) = [0,\:-1]\\
\end{array}
\]
where the infimum is w.r.t. to a switching sequence
$\sigma : [0,\infty) \mapsto \{1,\:2\}$ and
\[
A_1 := \left[\begin{array}{cc}-1&2\\1&-3\end{array}\right],
\quad A_2 := \left[\begin{array}{cc}-2&-2\\1&-1\end{array}\right].
\]
Since our framework cannot directly accomodate infinite-horizon
problems, we introduce a terminal condition $\|x(T)\|^2_2 \leq 10^{-6}$
so that terminal time $T$ is finite. It means that the switching
sequence should drive the state in a small ball around the origin.

\begin{table}[h!]
\centering
\begin{tabular}{c|c|c|cc}
$d$ & $p^*_d$ & $N_d$ & $y^d_{1,0}$ & $y^d_{2,0}$ \\ \hline
1 & $0.24294$ & 30 & 1.4252 & 1.4489 \\
2 & $0.24340$ & 105 & 2.0639 & 1.9237 \\
3 & $0.24347$ & 252 & 1.9639 & 1.8922 \\
4 & $0.24347$ & 495 & 1.9537 & 1.8904 \\
5 & $0.24347$ & 858 & 1.9572 & 1.8872 \\
6 & $0.24347$ & 1365 & 1.9677 & 1.8940  \\
7 & $0.24347$ & 2040 & 1.9669 & 1.8928  \\
\end{tabular}
\caption{Lower bounds $p^*_d$ on the optimal value $p^*$ obtained
by solving LMI relaxations of increasing orders $d$;
$N_d$ is the number of variables in the LMI problem;
$y^d_{k,0}$ is the approximate time spent on each mode $k=1,2$.\label{ex3table}}
\end{table}

In Table \ref{ex3table} we report the lower bounds $p^*_d$ on
the optimal value $p^*$ obtained by solving LMI relaxations of
increasing orders $d$, rounded to 5 significant digits.
We also indicate the number of variables
(i.e. total number of moments) of each LMI problem, as well as
the zeroth order moment of each occupation measure
(recall that these are approximations of the time spent
on each mode). We observe that the values of the
lower bounds stabilize quickly.

\begin{figure}[h!]
\centering
\includegraphics[width=\columnwidth]{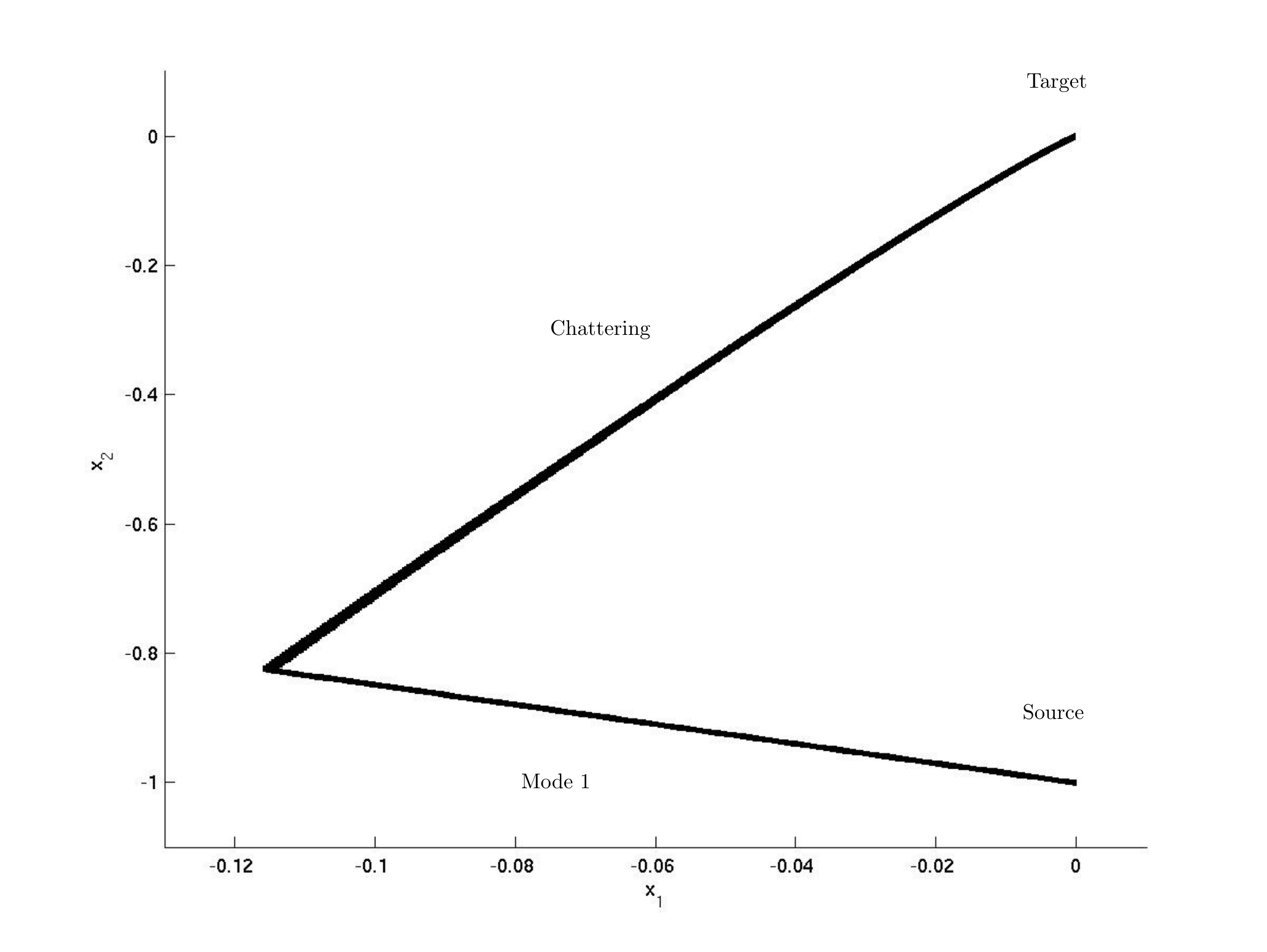}
\caption{Suboptimal trajectory starting at source point $x=(-1,0)$ with
mode 1, then chattering between modes 1 and 2 to reach the target,
a neighborhood of the origin.\label{ex3fig}}
\end{figure}

In Figure \ref{ex3fig} we plot an almost optimal trajectory
inferred from the moments of the occupation measure,
for an LMI relaxation of order $d=8$. The trajectory consists in starting
from $x=(-1,0)$ with mode 1 during 0.065 time units,
and then chattering between mode 1 and mode 2 with
respective proportions 49.3/50.7 until $x$ reaches
the neighborhood of the origin $\|x(T)\|^2_2 \leq 10^{-6}$
for $T=3.84$. This trajectory is slightly suboptimal,
as it yields a cost of $0.24351$, slightly bigger
than the guaranteed lower bound of $0.24347$ on
the best achievable cost obtained by the LMI relxation.
It follows that this trajectory is very close to optimality. 
For comparison with available suboptimal solutions, using
\cite[Theorem 1]{Geromel}, the so-called min switching strategy
yields with piecewise quadratic Lyapunov functions
a suboptimal trajectory with a cost of 0.24948.

\section{Conclusion}

In this paper we address the problem of designing an optimal switching sequence
for a hybrid system with polynomial Lagrangian (objective function)
and polynomial vector fields (dynamics). With the help of occupation measures,
we relax the problem from (control) functions with values in $\{0,1\}$
to (control) measures which are absolutely continuous w.r.t. time and
summing up to one. This allows for a convex linear programming (LP) formulation
of the optimal control problem that can be solved numerically
which a classical hierarchy of finite-dimensional convex
linear matrix inequality (LMI) relaxations.

We can think of two simple extensions of our approach:
\begin{itemize}
\item {Open-loop versus closed-loop}.
In problem (\ref{ocp}) the control signal is the switching sequence $\sigma(t)$ which is
a function of time: this is an open-loop control, similarly to what was proposed in \cite{impulse}
for impulsive control design. In addition, 
constrain the switching sequence to be an explicit or implicit function of the state,
i.e. $\sigma(x(t))$, a closed-loop control signal. In this case, each occupation measure
will be explicitly depending on time, state and control, and it will disintegrate as
$\mu(dt,dx,du)=\xi(dt\:|\:t,u)\omega(du\:|\:t)dt$, and we should follow the framework
described originally in \cite{sicon}.
\item {Switching and impulsive control}.
We may also combine switching control and impulsive control
if we extend the system dynamics (\ref{omega}) to
\[
dx(t) = \sum_{k=1}^m f_k(x(t)) \omega_k(dt) + \sum_{j=1}^{p} g_j(t) \tau_j(dt)
\]
where $g_j$ are given continuous vector functions of time
and $\tau_j$ are signed measures to be found, jointly with
the switching measures $\omega_k$. Whereas switching control measures $\omega_k$ are
restricted by (\ref{lebesgue}) to be absolutely continuous w.r.t.
the Lebesgue measure of time, impulsive control measure $\tau_j$ can
concentrate in time. For example, for a dynamical system
$dx(t) = g(t)\tau(dt)$, a Dirac measure $\tau(dt)=\delta_{s}$ 
enforces at time $t=s$ a state jump $x^+(s) = x^-(s) + g(s)$.
In this case, to avoid trivial solutions, the objective function
should penalize the total variation of the impulsive control measures,
see \cite{impulse}.
\item {Removing the states in the occupation measures}.
In the case that all dynamics $f_k(t,x)$, $k=1,\ldots,m$
are affine in $x$, we can numerically integrate the state
trajectory and approximate the arcs by polynomials
of time. It follows that the occupation measures $\mu_k(dt,dx)$,
once integrated, do not depend on $x$ anymore. They depend
on time $t$ only. We can then use finite-dimensional
LMI conditions which are necessary and sufficient for a vector to
contain the moments of a univariate measure, there is no need
to construct a hierarchy of LMI relaxations. There is however
still a hierarchy of LMI problems to be solved, now indexed
by the degree of the polynomial approximation of the arcs
of the state trajectory. To cope with high degree univariate
polynomials, alternative bases than monomials are recommended
(e.g. Chebyshev polynomials), see \cite{ltvimpulse} for more details.
\end{itemize}

\section*{Acknowledgments}

This work benefited from discussions with Milan Korda,
Jean-Bernard Lasserre and Luca Zaccarian.

\end{document}